\documentclass[11pt]{article}
\usepackage{amsmath,amsfonts,epsfig}
\newcommand{\bb}{\mathbb}

\newcommand{\cx}{{\bb C}}

\newcommand{\integers}{{\bb Z}}

\newcommand{\reals}{{\bb R}}

\newcommand{\hthree}{{\bb H}^3}

\newcommand{\ddz}{\frac{\partial}{\partial z}}
\newcommand{\qed}[1]{\nopagebreak[4]{\tiny \hfill
\fbox{\ref{#1}} \linebreak }\pagebreak[2]}
\newcommand{\arrow}{\rightarrow}
\newcommand{\del}{\partial}

\newcommand{\zbar}{{\overline{z}}}

\newcommand{\chat}{\widehat{\cx}}

\newcommand{\Ad}{\operatorname{Ad}}

\newcommand{\id}{\operatorname{id}}

\renewcommand{\Im}{\operatorname{Im}}

\newcommand{\ex}{\operatorname{ex}}

\newtheorem{theorem}{Theorem}[section]
\newtheorem{prop}[theorem]{Proposition}
\newtheorem{lemma}[theorem]{Lemma}
\newtheorem{cor}[theorem]{Corollary}

\newcommand{\cC}{{\cal C}}
\newcommand{\cD}{{\cal D}}

\newcommand{\cF}{{\cal F}}
\newcommand{\cG}{{\cal G}}

\newcommand{\cI}{{\cal I}}

\newcommand{\cL}{{\cal L}}

\newcommand{\cQ}{{\cal Q}}
\newcommand{\cR}{{\cal R}}

\newcommand{\dev}{D}
\newcommand{\torus}{E}

\newcommand{\upper}{U}

\newcommand{\bpq}{{\bf p},{\bf q}}
\begin{document}

\title{Hyperbolic Dehn surgery on geometrically infinite 3-manifolds}

\author{K. Bromberg\footnote{This work was partially supported by a
grant from the Rackham School of Graduate Studies, University of
Michigan and partially conducted for the Clay Mathematics Institute}}

\date{July 24th, 2000}

\maketitle

\begin{abstract}
In this paper we extend Thurston's hyperbolic Dehn surgery theorem to
a class of geometrically infinite hyperbolic 3-manifolds. As an
application we prove a modest density theorem for Kleinian groups. We
also discuss hyperbolic Dehn surgery on geometrically finite hypebolic
cone-manifolds.
\end{abstract}

\section{Introduction}
An essential result in the theory of hyperbolic 3-manifolds is
Thurston's Dehn surgery theorem. It states that all but a finite
number of Dehn fillings on a hyperbolic knot complement have a
hyperbolic structure. This theorem along with its extension to
geometrically finite manifolds has strong geometric and topological
consequences. In this paper we extend the Dehn surgery theorem to a
large class of geometrically infinite hyperbolic manifolds. As an
application we show that such geometrically infinite manifolds with
rank one cusps can be approximated by manifolds without rank one
cusps.

The Dehn filling theorem first appeared in Thurston's notes,
\cite{Thurston:book:GTTM}. A number of the pieces of the argument have
been published in a more rigorous form. We will indicate references as
we outline the result. There is also an approach to the finite volume
Dehn filling theorem via ideal triangulations. See
\cite{Petronio:Porti:dehn}.

In the simplest case, $M$ is a hyperbolic structure on the interior of
a compact manifold with boundary consisting of a single torus. This
hyperbolic structure defines a holonomy representation, $\rho$, of
$\pi_1(M)$ in $PSL_2\cx$. The first step is to show that the space of
all representations, $R(M)$, has complex dimension $\geq 1$ at $\rho$
(see \cite{Culler:Shalen:varieties}).

The complex length of an element in $PSL_2\cx$ measures the
translation length and rotation of the isometry. For an element in
$\pi_1(M)$ the complex length is a holomorphic function on $R(M)$. The
complex length of a curve on the torus boundary is zero at any
complete structure. Therefore by Mostow-Prasad rigidity this
holomorphic function is not constant and by the dimension count for
$R(M)$ must be locally onto.

To turn algebra into geometry, Thurston proves a
result that originated in work of Weil, \cite{Weil:compact}. Namely, if
$N$ is a compact hyperbolic manifold with holonomy $\rho$, and $\rho'$
is a representation near $\rho$, then there is a nearby hyperbolic
structure, $N'$, with holonomy $\rho'$. Furthermore, if $N''$ is
another nearby structure with holonomy $\rho'$, then there is a
slightly larger hyperbolic structure in which both $N'$ and $N''$
isometrically embed. A more careful statement and proof can be
found in \cite{Canary:Epstein:Green}.

For the Dehn filling theorem we let $N$ be a compact core of our
finite volume cusped manifold, $M$, and apply the above result. By the
algebraic result we know that the representation deforms and this in
turn gives a new hyperbolic structure, $N'$, on the compact core. The
above deformation theorem does not apply to the non-compact
cusp. Instead the geometry of the torus end must be constructed by
hand.

The metric completion of the structures we build on the torus ends is
described by a pair of real numbers, $(a,b)$, the {\em hyperbolic Dehn
filling coefficient}. If $a/b$ is rational
then the metric completion will consist of adding a simple close
curved. The topology of the completion will then be a solid torus. Let
$(p,q)$ be a relatively prime pair of integers such that $a/b =
p/q$. Then the $(p,q)$-curve on the torus will be trivial in the solid
torus. The metric around the added curve will usually be singular with
cross section a hyperbolic cone. If $p=a$ and $q=b$ then the metric
will be smooth.

Once the structures on the torus end have been constructed we need a
gluing theorem to show that the hyperbolic structure, $N'$, on the
compact core can be extended to one of these structures. This part of
Thurston's argument has not been published. Although not difficult,
the application of this result to geometrically infinite ends gave the
the author enough pause that he felt it should be written down in
detail. 

If we assign the complete, cusped structure the Dehn
filling coefficient $\infty$, then the map to the space of Dehn filling
coefficients, $\reals^2 \cup \infty$, is continuous. The final step in
the Dehn filling theorem is to show that if the complex length
functions for curves on the torus boundary is locally onto, then the
map to Dehn filling coefficients is locally onto. Since all but a
finite number of relatively prime pairs are in a neighborhood of
$\infty$ this completes our outline.

We now assume that $M$ is a complete hyperbolic structure on the
interior of a compact manifold with boundary consisting of two
components, a torus and  a surface of higher genus, $S$. We also
assume that this structure is geometrically finite without rank one
cusps. In this case $R(M)$ will have complex dimension $\geq 1 + \dim
T(S)$ where $T(S)$ is the Teichm\"uller space of $S$. The dimension of
the space of complete structures on $M$ is equal to the dimension of
$T(S)$, so once again the map to Dehn filling coefficients is locally
onto.

Every Dehn filling coefficient near $\infty$ is realized by a
representation $\rho'$ near $\rho$. Once again this gives a hyperbolic
structure on a compact core $N'$. This structure extends at the torus
end just as before. At the higher genus end we wish to extend $N'$ to
a geometrically finite structure. This is possible because
geometrically finite ends have a strong stability property: For any
small perturbation of the 
representation of a geometrically finite end, there is a nearby
geometrically finite end with that representation. Then Dehn filling
theorem for geometrically finite manifolds is proved in
\cite{Comar:thesis}.  Also see \cite{Bonahon:Otal:shortgs}.

In \cite{Hodgson:Kerckhoff:cone}, Hodgson and Kerckhoff extend the
Dehn filling theorem to cone-manifolds. In particular, they prove that
finite volume cone-manifolds are locally rigid when all cone angles
are $\leq 2\pi$. Furthermore, they show that local rigidity implies
that the map to Dehn filling coefficients is a local homeomorphism,
not just locally onto. This is extended to geometrically finite cone
manifolds in \cite{Bromberg:rigidity}. This algebraic result allows us
to prove the following strong Dehn filling theorem: Geometrically
finite cone-manifolds are locally parameterized by Dehn filling
coefficients and the conformal structure at infinity if all cone
angles are $\leq 2\pi$. The complete case corresponds to cone angle $0$.

The main focus of this paper is extending the Dehn filling theorem to
geometrically infinite 3-manifolds. If $M$ is geometrically infinite
many parts of the argument remain the same. The dimension count for
$R(M)$ will still 
hold and the map to Dehn filling coefficients will still be onto.
However, the geometrically infinite higher genus end will not be
stable. To get around this problem we restrict to a class of
structures where the geometrically infinite end will be
semi-stable. If none of the curves on the boundary torus of $M$ are
homotopic to the geometrically infinite end then there will be a half
dimensional subspace of  the representation variety for $S$,
where the infinite end will be stable. What we will see is that for
every Dehn filling coefficient near $\infty$ there will be a
representation that restricted to the geometrically infinite will be
in the stable subspace.

The main application of the geometrically infinite Dehn filling
theorem will be an approximation theorem for Kleinian groups. As an
example let $\Gamma$ be a Kleinian group such that $\hthree/\Gamma =
M$ is homeomorphic to the interior of a manifold with two
incompressible boundary components. We assume that one of the ends of
$M$ is geometrically infinite while the other is geometrically finite
{\em with} rank one cusps. Our approximation theorem will show that we
can resolve the rank one cusps. That is we can find Kleinian groups,
$\Gamma_i$, and hyperbolic 3-manifolds, $M_i = \hthree/\Gamma$, such
that $\Gamma_i \rightarrow \Gamma$ and $M_i$ is homeomorphic to $M$
with one geometrically infinite quasi-isometric to the infinite end of
$M$ and the other end geometrically finite {\em without} rank one
cusps. This is a very special case of the density conjecture of Bers
and Thurston.

In \S 2 we collect the background that we will need. In \S 3, we prove
the algebraic piece of the geometrically infinite Dehn filling
theorem. Then next two sections are expositions of material than can be
found in Thurston's notes: In \S 4 we prove the gluing theorem and in
\S 5 we construct the structures on the torus ends. For the Dehn
filling theorems for infinite volume hyperbolic manifolds we need to
construct structures on higher genus ends. This is done in \S 6. We put
everything together to prove the Dehn filling theorems in \S7. We prove
the approximation theorem in \S 8.

{\bf Acknowledgments.} The author would like to thank Dick Canary
for many helpful discussions.

\section{Preliminaries}

We begin with a review of some basic facts about hyperbolic 3-manifolds.

\subsection{Kleinian groups}
Hyperbolic 3-space, $\hthree$, is the unique simply connected
3-manifold with constant sectional curvature equal to $-1$. The group of
orientation preserving isometry groups of $\hthree$ is naturally
isomorphic to $PSL_2\cx$. Furthermore $\hthree$ can be compactified
by the Riemann sphere, $\chat$, and the isometric action of elements
of $PSL_2\cx$ on $\hthree$ extends continuously to projective
transformations of $\chat$.

In this paper,
a {\em Kleinian group}, $\Gamma$, is a discrete, torsion free subgroup of
$PSL_2\cx$.  $\Gamma$ will act properly
discontinuously on $\hthree$ so the quotient, $M = \hthree/\Gamma$, will be a
hyperbolic 3-manifold.

\subsection{The domain of discontinuity and limit set}
The action of $\Gamma$ on $\chat$ naturally decomposes $\chat$ into
two sets. The {\em domain of discontinuity}, $\Omega(\Gamma)$, is the
largest subset of $\chat$ where $\Gamma$ acts properly
discontinuously. The {\em limit set}, $\Lambda(\Gamma)$, is the
complement of $\Omega(\Gamma)$. By Ahlfors finiteness theorem, if
$\Gamma$ is finitely generated,
$\Omega(\Gamma)/\Gamma$ is a finite collection of Reimann surfaces of
finite type. As we shall see $\Omega(\Gamma)/\Gamma$ can be naturally
identified with the quotient hyperbolic 3-manifold to form a (not
necessarily compact) 3-manifold with boundary.

\subsection{Laminations}
There are a number of ways to define a lamination on a surface
$S$. The easiest is to endow $S$ with a hyperbolic metric. Then a
{\em lamination} on $S$ is a collection of disjoint, simple geodesics whose
union is a closed subset of $S$. However, the use of the hyperbolic
metric is unnecessary and with more work we could define this concept
pure topologically.

A lamination, $\lambda$, is finite if it consists entirely of simple
closed curves on $S$. In general, $\lambda$, will be a limit of finite
laminations in the Hausdorff topology on closed subsets of $S$.

\subsection{Ends} Throughout this paper we will assume that all 3-manifolds
are tame. That is, $M$ will be the interior of a compact 3-manifold
with boundary, $N$. If $S$ is a connected component of $\del N$ then a
neighborhood $E$ of $S$ is an {\em end} of $N$. In abuse of notation $E$
will also refer to the restriction of this neighborhood to $M$. We can
assume that $E$ is an $I$-bundle over $S$. In particular we
identify $E$ with $S \times [0,1]$ with $S \times 1$ lying in $\del N$.

A sequence of simple closed curves, $c_n$, {\em
exits} $E$ if for every $t \in [0,1)$, $c_n$ is contained in $S \times
[t,1)$ for all large $n$.

\subsection{Ending laminations} A complete hyperbolic structure on $M$
determines a lamination, $\lambda$, on each component $S$ of $\del
N$. Loosely speaking $\lambda$ will be the collection of geodesics
which exit the end $E$. We make a precise definition of $\lambda$ by
examining its connected components.

If a connected component, $c$, of $\lambda$ is a simple closed
curve then there will be simple closed curves $c_i$ exiting $E$ each
isotopic to $c$ in $N$ such that the length of the $c_i$ in $M$ goes
to zero.

For a connected component, $\lambda_0$, of $\lambda$ that is not a
simple closed curve there will be a sequence of simple closed curves
$c_i$ on $S$ converging to $\lambda_0$ such that the geodesic
representatives of the $c_i$ exit $E$. We say that $c$ is a {\em
rank-one cusp} of $E$.

The union of all such laminations satisfying these two conditions is
the {\em ending lamination}, $\lambda$, for $E$. For a tame end the work of
Bonahon and Canary shows that $\lambda$ is well defined and
unique. The ending lamination for $M$ is the union of the ending
laminations for each end.

\subsection{The conformal structure at infinity} The Kleinian manifold
$\hat{M}$ is the quotient of $\hthree \cup \Omega(\Gamma)$ under the
action of $\Gamma$. The inclusion of $M$ into $N$ will extend to an
inclusion of $\hat{M}$ into $N$. The image of $\hat{M}$ in $\del N$
will be a collection of essential subsurfaces bounded by simple closed
curves in the ending lamination. Since $\Omega(\Gamma)/\Gamma$ is a
Riemann surface we have defined a {\em conformal structure at
infinity} on a subsurface of $\del N$.

\subsection{End-invariants} We have now defined two objects on $\del N$: the
ending lamination, $\lambda$, and the conformal structure at infinity,
$\mu$. The key fact is that these two objects will have disjoint
support. Furthermore $\lambda$ will be a maximal in the sense that any
simple closed curve disjoint from $\lambda$ and the support of $\mu$
will be homotopically trivial. The {\em end-invariant} of a complete
hyperbolic structure on $M$ is the pair $(\lambda, \mu)$. Thurston's
ending lamination conjecture states that this pair uniquely determines
the hyperbolic structure on $M$.

\subsection{Accidental parabolics} If $E$ is the end associated to the
component $S$ of $\del N$ then $c$ is an {\em accidental parabolic} if
$c$ has parabolic holonomy yet is not in the ending lamination of
$E$. As a consequence of negative curvature, two simple closed curves
that are not homotopic on $\del N$ but are homotopic in $N$ cannot
both be in the ending lamination. Therefore if one of the curves is in
the ending lamination the other will be an accidental parabolic.

\subsection{Geometrically finite and infinite hyperbolic structures} 
\label{findef}
If an end has empty ending lamination then it is {\em geometrically
finite without rank-one cusps}. If the ending lamination is a finite
lamination then the end is simply {\em geometrically
finite}. Otherwise the end is {\em geometrically infinite}. We make
similar definitions for hyperbolic structures on the entire manifold,
$M$.

\subsection{Quasiconformal deformation spaces}
A {\em quasiconformal homeomorphism} between two Riemann surfaces is a
homeomorphism, $f$, with distributional derivatives locally in $L^2$
and $\|f_\zbar/f_z\|_\infty < 1$. The quantity, $\mu = f_\zbar/f_z$ is
the {\em beltrami differential} for $f$ and is a differential of type
$(-1,1)$. The map $f$ is $k$-quasiconformal if $k =
\frac{1+\|\mu\|_\infty}{1 - \|\mu\|_\infty}$.

Let $\Gamma$ and $\Gamma'$ be two Kleinian groups with a group
isomorphism taking $\gamma \in \Gamma$ to $\gamma' \in
\Gamma'$. $\Gamma'$ is a {\em quasiconformal deformation} of $\Gamma$
if there is a quasiconformal homeomorphism, $f: \chat \longrightarrow
\chat$, such that $f \circ \gamma = \gamma' \circ f$ for all $\gamma
\in \Gamma$.
Let $\cQ\cD(\Gamma)$ be the space of Kleinian groups which are
quasiconformal deformations of $\Gamma$ and $QD(\Gamma)$ be the
quotient of $\cQ\cD(\Gamma)$ under the action of $PSL_2\cx$ by
conjugation.

Assume $M = \hthree/\Gamma$ is tame and let $S$ be a subsurface of
$\del M$ bounded by simple closed curves in the ending
lamination. Then a component of $\Omega(\Gamma)/\Gamma$ lie in $S$ in
we let $\Omega_S \subseteq \Omega(\Gamma)$ be the pre-image of this
component. Then $QD(\Gamma;S) \subseteq QD(\Gamma)$ is the space of
quasiconformal deformations such that there exist a conjugation
quasiconformal homeomorphism whose beltrami differential's support is
disjoint from $\Omega_S$.

It will often be convenient to refer to the elements of $QD(\Gamma)$
which are equivalence classes of groups by a single group in the
class. We define $\Gamma' \in QD(\Gamma)$ to the the unique group in
it's equivalence class that is the conjugation of $\Gamma$ by a
quasiconformal homeomorphism fixing $0$, $1$ and $\infty$. This also
gives a canonical identification of $\cQ\cD(\Gamma)$ with $QD(\Gamma)
\times PSL_2\cx$.

\subsection{Teichm\"uller spaces}
\label{teichspace}
Let $S$ be a closed surface with a finite
number of punctures. The {\em Teichm\"uller space} of $S$, $T(S)$, is
the space marked, complex structures on $S$ or equivalently the space
of marked, complete, finite area hyperbolic structures on $S$. $T(S)$
has a canonical complex structure and is a cell of complex dimension
$3g - 3 + p$ where $g$ is the genus of $S$ and $p$ the number of
punctures.

The following theorem is due to Ahlfors, Bers, Kra, Maskit and Sullivan.
\begin{theorem}
\label{locpar}
$QD(\Gamma)$ is locally parameterized by $T(\Omega(\Gamma)/\Gamma)$.
$QD(\Gamma;S)$ is locally parameterized by $T(\Omega_C/\Gamma)$ where
$\Omega_C$ is the complement of $\Omega_S$ in $\Omega(\Gamma)$.
\end{theorem}
This parameterization is global if $\hthree/\Gamma$ has incompressible
boundary. There is also a global parameterization for the general case
although it is more difficult to state and we will not need it.

Note that the complex structure on Teichm\"uller space defines a complex
structure on $QD(\Gamma)$.

\subsection{Representation varieties}
We let $\cR(\Gamma)$ be the space all representations of the abstract
group $\Gamma$ in $PSL_2\cx$. Since $PSL_2\cx$ is a complex, algebraic
group, $\cR(\Gamma)$ has a naturally structure as a complex
variety. We let $R(\Gamma)$ be the quotient of $\cR(\Gamma)$ under
the action of $PSL_2\cx$. If $\Gamma$ is the fundamental group of a
manifold $M$, $\cR(M)$ and $R(M)$ will be representations
of $\pi_1(M)$ and the quotient space, respectively.

There is a natural inclusions of $QD(\Gamma)$ in $R(\Gamma)$. The
Ahlfors-Bers version of the measurable Reimann mapping theorem implies
that this inclusion is holomorphic. Kapovich 
has shown that $R(\Gamma)$ is a smooth variety at all discrete
faithful representations. In particular, $R(\Gamma)$ is a smooth
complex manifold in a neighborhood of the image of
$QD(\Gamma)$. Therefore we can view the inclusion map as a
holomorphic embedding of one complex manifold in another.

\subsection{Group cohomology}
A {\em 1-cocyle} is a map, $z:\Gamma \longrightarrow sl_2\cx$ with
$z(\gamma_1\gamma_2) = z(\gamma_1) + \Ad \rho(\gamma_1) z(\gamma_2)$. A
1-cocycle is a {\em 1-coboundary} if there exists a $v \in sl_2\cx$
such that $z(\gamma) = v - \Ad \rho(\gamma) v$ for all $\gamma \in
\Gamma$. We let $Z^1(\Gamma; \Ad \rho)$ be the space of 1-cocycles,
$B^1(\Gamma; \Ad \rho)$ the space of 1-coboundaries, and $$H^1(\Gamma;
\Ad \rho) = Z^1(\Gamma; \Ad \rho)/B^1(\Gamma; \Ad \rho).$$ We then
have: 
\begin{theorem}[\cite{Weil:cohomology}]
\label{weiltangent}
$H^1(\Gamma; \Ad \rho)$ is canonically isomorphic to the Za\-ri\-ski
tangent space, $TR(\Gamma)$, of  $R(\Gamma)$ at $\Gamma$.
\end{theorem}

{\bf Remark.} Strictly speaking, for this theorem to be true,
$H^1(\Gamma; \Ad \rho)$ should be viewed as the Zariski tangent space
of the scheme, $R(\Gamma)$. However, at the discrete faithful
representations, Kapovich's theorem, mentioned above, implies that
$R(\Gamma)$ is a smooth complex manifold and the Zariski tangent space
can be canonically identified with the differentiable tangent space.

\subsection{Strain fields}
For a vector field on $\chat$ the {\em strain} is the infinitesimal
version of the beltrami differential of a homeomorphism. As we will also
need the notion of the strain of a vector field on $\hthree$, we
define the strain, $Sv$, of a general vector field, $v$, on a
conformally flat manifold as a tensor of type (1,1) with
$$S_{ij}v =
\frac12\left(\frac{\del v_i}{\del x_j} + \frac{\del v_j}{\del
x_i}\right) - \frac{\delta^i_j}{n}\Sigma^n_{i=1} \frac{\del v_i}{\del
x_i}$$
where $v = (v_1, \dots, v_n)$ on a local, conformal chart. If $v =
f\ddz$ is a vector field on $\chat$ an easy calculation shows that $Sv
= f_\zbar \frac{d\zbar}{dz}$. We also let $|Sv|$ be the operator norm
of $Sv$. Then $v$ is a {\em quasiconformal vector field} if $Sv$ is
measurable and $|Sv|$ is essentially bounded. 

\subsection{Infinitesimal deformations}
An {\em infinitesimal quasiconformal deformation} of $\Gamma$ is a
quasiconformal vector field, $v$, on $\chat$ such that there exists a
map, $z:\Gamma \longrightarrow sl_2\cx$ with $v - \gamma_* v=
z(\gamma)$. By the chain rule $z$ is cocycle.
An {\em infinitesimal hyperbolic deformation} of $\Gamma$ is a
smooth vector field, $v$, on $\hthree$ such that $v - \gamma_*v =
z(\gamma)$. 

Note that both types of infinitesimal deformations form a vector space
and there is a homomorphism from each of these vector spaces to
$H^1(\Gamma; \Ad \rho)$. For infinitesimal hyperbolic deformations
this map will be onto. In particular, we can define $v$ to be zero on
a fundamental domain for $\Gamma$ in $\hthree$ and use the cocycle
condition to define $v$ on all of $\hthree$. The vector field, $v$,
can be made smooth via a partition of unity argument.

On the other hand the map from infinitesimal quasiconformal
deformations to $H^1(\Gamma; \Ad \rho)$ will not always be onto. Let
$TQD(\Gamma)$ be the 
subspace of $TR(\Gamma)$ which is the tangent space of the embedding
of $QD(\Gamma)$ in $R(\Gamma)$. We then have the following theorem
which can be found in \cite{Ahlfors:finiteness}.
\begin{theorem}
\label{qctangspace}
Let $z \in TQD(\Gamma)$ be a cocycle. Then there exists an
infinitesimal quasiconformal deformation, $v$, such that $v - \gamma_*
v = z(\gamma)$. 
\end{theorem}

\subsection{Continuous extension of vector fields on $\hthree$ to $\chat$}
If we view $\hthree$ as the open unit ball in $\reals^3$ and $\chat$
as the unit sphere then we can discuss the continuous extension of
vector fields on $\hthree$ to $\chat$. A general vector field on
$\hthree$ will not continuously extend to $\chat$, and even if the
vector field in the unit ball does extend to the unit sphere, it may
not be tangent to the unit sphere and therefore will not be a
continuous vector field on $\chat$. This leads us to consider the
weaker notion of a continuous {\em tangential} extension.  Given a
vector field $v$ on $\hthree$ and a vector field $V$ on $\chat$ we say
that $V$ is a continuous tangential extension of $v$ if there exists a
point $p \in \hthree$ such that the projection of $v$ to the family of
hyperbolic spheres centered at $p$ extends continuously to
$\chat$. Note that the choice of $p$ is arbitrary and if such a
continuous extension exists for a single point, $p \in \hthree$, then
it will exist for every point in $\hthree$ and the extension to $\chat$
will always be the same vector field.

The following lemma is an easy exercise:
\begin{lemma}
\label{autext}
If $v$ is an infinitesimal hyperbolic deformation and $V$ is a
continuous tangential extension of $v$ that is a quasiconformal vector
field then $V$ will be an infinitesimal quasiconformal deformation
that determines the same cocycle as $v$.
\end{lemma}

\subsection{Averaging vector fields on $\chat$}
Let $V(\chat)$ and $V(\hthree)$ be the space of (distributional)
vector fields on 
$\chat$ and $\hthree$ respectively. The following theorem is due to
Ahlfors, Reimann and Thurston. A self-contained exposition can be
found in \cite{McMullen:book:RTM}.
\begin{theorem}
\label{average}
There exists a unique operator, $\ex: V(\chat) \longrightarrow
V(\hthree)$, such that for every $v \in V(\chat)$ the following
properties hold:
\begin{enumerate}
\item $\ex(v)$ is smooth and extends continuously to $v$ if $v$ is
continuous.

\item If $v$ is quasiconformal, $\ex(v)$ is quasiconformal.

\item For every $\gamma \in PSL_2\cx$, $\gamma_*\ex(v) =
\ex(\gamma_*v)$.

\item If $v$ is a projective vector field, $\ex(v)$ is an
infinitesimal isometry.
\end{enumerate}
\end{theorem}
The vector field $\ex(v)$ evaluated at $x$ is essentially the the {\em
visual average} of $v$ when viewed from $x$ in hyperbolic space.

\section{Geometrically infinite manifolds}
\label{infinite}
Let $N$ be a compact, hyperbolizable 3-manifold with interior
$M$. Assume $N$ is bounded by $k$ tori and $n-k$ surfaces of genus
$> 1$ with $n>k>1$. Label the tori $T_i,
\dots, T_k$ and the surfaces of higher genus $S_{k+1}, \dots, S_n$ and
let $E_i$ be the associated ends. Let $\Gamma$ be a Kleinian group
such that $M = \hthree/\Gamma$.  Let $\tilde{E}_i$ be a lift
of $E_i$ to the universal cover, $\tilde{M}$, and $\Gamma_i \subset
\Gamma$ the subgroup fixing $\tilde{E}_i$. Choose an integer $m$, such
that $k<m\leq n$. Throughout this section we make the following
assumptions:

\begin{enumerate}
\item The surfaces $S_i$ are incompressible for $k< i \leq m$.

\item $E_i$ does not have accidental parabolics for $k <i \leq m$.

\item The hyperbolic structures on $E_i$ are geometrically finite
without rank one cusps for $i > m$.

\end{enumerate}

Conditions (1) and (2) restrict the topological type of the manifold
$M$ along with the geometry. Although the class of manifolds that
satisfy these two conditions is quite large, some of the simplest
hyperbolic manifolds are not in it, such as $I$-bundles over a surface
with a curve removed. The importance of these conditions becomes
apparent in the following proposition.

\begin{prop}
\label{bgroup}
The dimension of $QD(\Gamma_i; S_i)$ is equal to the dimension of the
Teichm\"uller space of $S_i$ for $i \leq m$.
\end{prop}

{\bf Proof.} Since $S_i$ is incompressible, Bonahon's theorem
\cite{Bonahon:tame} implies
that $\hthree/\Gamma_i$ is homeomorphic to $S_i \times (0,1)$. On the
other hand, $M$ is not an I-bundle because it has two non-homeomorphic
boundary components. Therefore $\Gamma_i$ has infinite index in
$\Gamma$. We can assume $E_i$ has been chosen such that
the covering map is injective onto $E_i$. Let $E^*_i$ be the other end
of $\hthree/\Gamma_i$. The covering map from $E^*_i$ to $N$ must be
infinite-to-one and therefore by the covering theorem
\cite{Canary:inj:radius}, $E^*_i$ is geometrically finite.

Since $E^*_i$ is geometrically finite its ending lamination consists
entirely of simple closed curves. However, if $c$ is a simple
closed curve in the ending lamination for $E^*_i$ then $c$ is
parabolic. However by assumption (3), $c$ is then in the ending
lamination for $E_i$. This is a contradiction so $E^*_i$ has an empty
ending lamination, and by Theorem \ref{locpar}, $QD(\Gamma_i;S_i)$ is
locally parameterized by the Teichm\"uller space of $S_i$.
\qed{bgroup}

To simplify notation we let $\del_1 M = S_{k+1} \cup \cdots \cup S_m$ and
$\del_2 M = S_{m+1} \cup \cdots \cup S_{n}$. $R(\del_1 M)$, $R(\del_2 M)$ and
$T(\del_2 M)$ are then the products of the respective representation
varieties and Teichm\"uller spaces. $QD(\del_1 M)$ is the product of
relative quasiconformal deformation spaces $QD(\Gamma_i; S_i)$ with
$i=k+1,\dots,m$. Let $V$ be a small neighborhood $\rho$ in $R(\Gamma)$.
We then have restriction maps $\del_i:R(\Gamma) \longrightarrow R(\del_i
M)$ and for the geometrically finite ends a map to Teichm\"uller
space, $b'_2:V \longrightarrow T(\del_2 M)$ (see \S \ref{ends}).
For each torus boundary component a choice of a simple closed curve
determines a length function, $\cL'_i:V \longrightarrow \cx^k$
(see \S \ref{torus}).
Putting these maps together we have maps, $\del:R(\Gamma)
\longrightarrow \cx^k \times R(\del_1 M) \times R(\del_2 M)$ defined
by $\del = ( \cL_1,\dots,\cL_k,\del_1, \del_2)$ and $\Phi':V
\longrightarrow \cx^k \times T(\del_2 M)$ defined by $\Phi' = (\cL'_1,
\dots, \cL'_k, b'_2)$.

Finally we define $CD(\del_2 M) = (b'_2)^{-1}(b'_2(\del_2(\Gamma)))$, the
space of conformal deformations of the geometrically finite ends and
$R(\Gamma; \del_1 M) = \del^{-1}(\cx^k \times QD(\del_1 M) \times
R(\del_2 M))$. 

\begin{theorem}
\label{infinitedehn}
$\Phi'$ restricted to $R(M; \del_1 M)$ is a local homeomorphism
at $\Gamma$.
\end{theorem}

{\bf Proof.} We begin by noting the smoothness and dimension of the
various spaces. All dimensions are complex. $R(\Gamma)$ is smooth at
$\Gamma$ and has dimension $-\frac32 \chi(\del M) + k$ by Theorem
9.8.1 in \cite{Kapovich:book}. $R(\del_1 M)$ and $R(\del_2 M)$ are
smooth at $\del_1(\Gamma)$ and $\del_2(\Gamma)$, respectively, with
dimension $-3 \chi(\del_1 M)$ and $-3 \chi(\del_2 M)$, respectively
\cite{Gunning:book:LVB}. $QD(\del_1 M)$ is a complex submanifold of
$R(\del_1 M)$ and has dimension $-\frac32 \chi(\del_1 M)$ by Lemma
\ref{bgroup}. $CD(\del_2 M)$ is isomorphic to the space of holomorphic
quadratic differentials on $\del_2 M$ and therefore is smooth of
dimension $-\frac32 \chi(\del_2 M)$. Let $TR(\Gamma)$, $TR(\del_1 M)$,
$TR(\del_2 M)$, $TQD$, and $TCD$  be the respective tangent spaces at
$\Gamma$ or its image under the appropriate map.

To prove the theorem we need to show that $$\dim R(M; \del_1 M)
= \dim T(\del_2 M) + k = -\frac32 \chi(\del_2 M) + k$$ and that
$(\Phi')_*$ is injective.
%We first show that the second condition implies the first.

Assume that $v \in TR(\Gamma)$ is an infinitesimal hyperbolic
deformation. We make the following claim: If $\del_*v \in
(0 \times TQD \times TCD)$ then $v$ is a trivial deformation. We first
will assume this claim and complete the proof of the theorem, saving
the proof of the claim for below.

We first show that $R(M; \del_1 M)$ is smooth and has the correct
dimension by showing
that $\del_*TR(M)$ is transverse to $\cx^k \times TQD \times TR(\del
M_2)$ and that the dimension, $D$, of 
$$\del_*TR(\Gamma) \cap (\cx^k \times TQD \times TR(\del_2 M))$$
is $-\frac32 \chi(\del M_2) + k$. The claim implies that $\del_*$ is
injective and therefore 
$$\dim \del_* TR(\Gamma) = -\frac32 \chi(\del M) + k$$
while
\begin{eqnarray*}
\dim(\cx^k \times TQD \times TR(\del_2 M)) &=& -\frac32 \chi(\del_1
M) - 3\chi(\del_2 M) + k \\
&=& -\frac{3}{2}\chi(\del M) -\frac{3}{2}\chi(\del_2 M) + k
\end{eqnarray*}
and
\begin{eqnarray*}
\dim (\cx^k \times TR(\del_1 M) \times TR(\del_2 M))
& = & -3\chi(\del_1 M) - 3\chi(\del_2 M) + k \\
& = & -3\chi(\del M) + k
\end{eqnarray*}
so
$$D \geq -\frac32 \chi(\del_2 M) + k.$$
On the other hand, the claim also tells us that
$$\del_*TR(M) \cap (0 \times TQD \times TCD) = 0$$
and since
$$\dim (\cx^k \times TQD \times TCD) = -\frac32 \chi(\del_1 M) -
\frac32 \chi(\del_2 M) + k,$$
we have
$$D \leq -\frac32 \chi(\del_2 M) + k.$$
Therefore
$$D= -\frac32 \chi(\del_2 M) + k.$$
Since the dimension of intersection is minimal, the intersection is
transverse and therefore $R(\Gamma; \del_1 M)$ is smooth at $\Gamma$
and has dimension $D = \dim T(\del_2 M) + k$.

We now show that $(\Phi')_*$ restricted to $TR(\Gamma; \del_1 M)$ is
injective. Since $R(\Gamma; \del_1 M)$ is smooth we can discuss its
tangent space, $TR(\Gamma; \del_1 M)$, at $\Gamma$.  If $v \in TR(\Gamma;
\del_1 M)$ then $\del_*v \in (\cx^k \times TQD \times TR(\del_2 M))$
while if $(\Phi')_*v = 0$ then $\del_*v \in (\cx^k \times TQD \times
TCD)$. By the claim, $v$ is trivial, so the restriction of
$(\Phi')_*$ to $TR(\Gamma; \del_1 M)$ is injective and since
$TR(\Gamma; \del_1 M)$ and $\cx^k 
\times T(\del_2 M)$ have the same dimension, $(\Phi')_*$ is an
isomorphism and therefore $\Phi'$ restricted to $R(\Gamma; \del_1 M)$
is a local homeomorphism at $\Gamma$.

We now prove the claim in a sequence of two lemmas.

\begin{lemma}
\label{discdomainext}
If $\del_*v \in (0 \times TQD \times TCD)$ then $v$ can be represented
by a quasiconformal vector field that extends continuously to
$\Omega(\Gamma)$.
\end{lemma}

{\bf Proof.} We first construct the vector field on the lift of the
geometrically infinite ends. By assumption, $(\del_1)_*v \in TQD$ so
Theorem \ref{qctangspace} implies that there exists a
vector field $V_i$ on $\chat$ which is a infinitesimal quasiconformal
deformation equivalent to $v$. We then let $v' = \ex(V_i)$ on
$\tilde{E}_i$. Theorem \ref{average} implies that $v'$ is
quasiconformal on $\tilde{E}_i$. For $\gamma \in
\Gamma$, define $v' =  \gamma_* v' + v - \gamma_* v$. It is easy to
check that this is well defined using the fact that $v$ is an
infinitesimal deformation of $\Gamma$ and $v'$ is an infinitesimal
deformation of $\Gamma'$. Furthermore, where it is defined, $v'$ will
be an infinitesimal hyperbolic deformation defining the same cocycle
as $v$.

If $i > m$ we construct $v'$ on $\tilde{E}_i$ using Theorem 3.9 of
\cite{Bromberg:rigidity}. For the ends associated to the torus
boundary components we use Proposition 3.10 of
\cite{Bromberg:rigidity}.

The vector field, $v'$, is now defined outside of the universal cover
of a compact subset 
of $M$, so as long as the extension of $v'$ to the rest of $\hthree$ is
smooth, $v'$ will be quasiconformal. This extension can be done by
combining $v$ with $v'$ using a cutoff function. We again
refer to \cite{Bromberg:rigidity} for details. \qed{discdomainext} 

In the next lemma we show that this vector field extends to all of
$\chat$.

\begin{lemma}
\label{limitext}
Every quasiconformal vector field, $v \in TR(M)$, which extends
continuously to $\Omega(\Gamma)$ has a continuous tangential extension
to a quasiconformal vector field on all of $\chat$.
\end{lemma}

{\bf Proof.}
If $\Lambda(\Gamma) = \chat$ then $v$ has a continuous tangential
extension on $\chat$ by Theorem 4.8 in
\cite{Kapovich:extension}.  If not
we must	modify Kapovich's result only slightly.

Kapovich's proof has two steps:

1. The infinitesimal deformation, $v$, extends continuously to the
loxodromic fixed points of $\Gamma$.

2. If a quasiconformal vector extends continuously to a dense subset
of $\chat$ then it has a continuous tangential extension to a
quasiconformal vector field on all of $\chat$.

If $\Lambda(\Gamma) = \chat$ then we follow Kapovich and apply
(1). Since the loxodromic fixed points are dense in $\Lambda(\Gamma)$
this gives a continuous extension to a dense subset of $\chat$.
Otherwise $\Omega(\Gamma)$ is dense in $\chat$ and we have assumed that
$v$ extends continuously to $\Omega(\Gamma)$. In both cases we then
apply (2).
\qed{limitext}

It is now an easy matter to complete the proof of the claim. By Lemmas
\ref{discdomainext} and \ref{limitext} 
there exists a deformation, $v'$, equivalent to $v$ that has a
continuous tangential extension to a 
quasiconformal vector field, $V$, on $\chat$ that is conformal on
$\Omega(\Gamma)$. Furthermore, by Lemma \ref{autext}, $V$ is an
infinitesimal quasiconformal deformation that defines the same
cocycle as $v$. Therefore $SV$ will be an equivariant strain field
which is zero on $\Omega(\Gamma)$ and by Sullivan's rigidity theorem,
\cite{Sullivan:linefield},  $SV$ is 
zero on all $\chat$. This implies that $V$ is a conformal vector field
and that the associated cocycle is a coboundary. Hence, $v$ is a
trivial deformation proving the claim and Theorem
\ref{infinitedehn}. \qed{infinitedehn}

\section{Developing maps}

A {\em geometry} is a pair made up of a simply connected manifold,
$X$, with a real analytic, transitive group action $G$. The
``geometry'' of $X$ may 
be a Riemannian metric with $G$ the group of isometries. However, this
is not the only possibility. For example $G$ may be the group of
projective or affine transformations of $X$.

For a manifold, $M$, a $(G,X)$-{\em developing map}, $\dev$, is a local
homeomorphism from 
$\tilde{M}$ to $X$ such that there exists a {\em holonomy representation},
$\rho: \pi_1(M) \longrightarrow G$ with
$$ \dev(\gamma(x)) = \rho(\gamma)(\dev(x))$$
for all $\gamma \in \pi_1(M)$ and $x \in \tilde{M}$. This defines a
{\em $(G,X)$-structure} on $M$.
We will denote the space of developing maps for $M$ with the
$C^\infty$ topology by $\cD(M)$. 

$M$ is a {\em thickening} of a compact manifold $N$ if $M - N$ is
homeomorphic to $\del N \times (0,1]$. The following is Theorem 1.7.1
in \cite{Canary:Epstein:Green}. 

\begin{theorem} 
\label{weils theorem}
Let $M$ be a thickening of $N$ and $D_0$ a developing
map for $M$ with holonomy representation $\rho_0$.

\begin{enumerate}
\item There exists a continuous map, $D$, from a neighborhood, $V$, of
$\rho_0$ in $\cR(M)$ to $\cD(M)$ such that the holonomy of $D(\rho)$
is $\rho$.

\item For any $V$ as in 1, there is a neighborhood of $D_0$ in
$\cD(N)$ that is homeomorphic to a product $\cI \times V$, where $\cI$
is a neighborhood of the inclusion $N \hookrightarrow M$ in the space
of locally flat embeddings. In particular, if $(\iota, \rho) \in \cI
\times V$ and $\tilde{\iota}$ is the lift of $\iota$ to the universal
cover, then $(\iota, \rho) \mapsto D(\rho) \circ \tilde{\iota}$ under
this homeomorphism.
\end{enumerate}
\end{theorem}

Let $M$ be an open manifold and $N_0$, $N_1$, $N_2$ and $N_3$
submanifolds of $M$ such that $N_{i+1}$ is a thickening of $N_i$ for
$i=0,1,2$ and $M$ is the interior of a thickening of $N_3$. Then {\em
ends} of $M$ will be the finite number of components of $M \backslash
N_1$. Let $E$ be an end of $M$. Finally, we define two more
submanifolds of $M$; $C_1$ is the closure of $(N_2 \backslash N_1)
\cap E$ and $C_2$ is the closure of $(N_3 \backslash N_0) \cap E$.

\begin{theorem}
\label{weilopen}
Let $D$ be a developing map for $M$ and $D_3$ and $D_E$ the
restriction of $D$ to $N_3$ and $E$, respectively. Assume that
$D'_3$ and $D'_E$ are small deformations of $D_3$ and $D_E$,
respectively, such their restrictions to $C_1$ have the same
holonomy. Then there exists a developing map, $D'$, for $M$ such that
$D'|_{N_0} = D'_3|_{N_0}$ and $D'|_E = D'_E$.
\end{theorem}

{\bf Proof.} 
\begin{figure}[ht]
\begin{center}
\epsfig{file=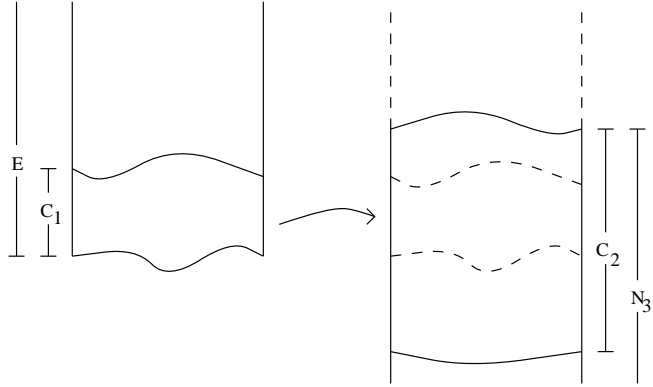,angle=-90}
\caption{After we deform $E$ and $N_3$, the new structure on $C_1$
coming from $E$ will embed in the new structure on $C_2$ coming from
$N_3$.}
\end{center}
\end{figure}
We apply Theorem \ref{weils theorem} to $C_1$ and
$C_2$. Let $\rho$ be the holonomy of $D$ restricted to $C_2$. Then, as
in (1), we have a map $D_C: V \longrightarrow \cD(C_2)$ where $V$ is a
neighborhood of $\rho$ in $\cR(C_2)$ and as in (2), we have a neighborhood
of $D$ restricted to $C_1$ parameterized by $V \times \cI$ where $\cI$
is a neighborhood of the inclusion map in the space of locally flat
embeddings. 

Let $\rho'$ be the holonomy of $D_3$ and $D_E$ restricted to $C_1$.
By (2) there exist embeddings, $\iota_3, \iota_E:C_1 \longrightarrow
C_2$, such that $D'_3|_{C_1} = D_C(\rho') \circ \tilde{\iota}_3$ and
$D'_E|_{C_1} = D_C(\rho') \circ \tilde{\iota}_E$. Extend
$\iota_3$ to a diffeomorphism, $\psi: N_3 \longrightarrow N_3$ such
that $\psi|_{N_0}$ is the identity. We can then find another
diffeomorphism, $\phi:N_3 \longrightarrow N_3$, such that $(\psi \circ
\phi) |_{C_1} = \iota_E$. Define the developing map, $D'$, such that
$D'$ restricted to $N_2$ is $D'_3 \circ \tilde{\psi} \circ
\tilde{\phi}$ and $D'$ restricted to $E$ is $D'_E$. Since $D'_E = D'_3
\circ \tilde{\psi} \circ \tilde{\phi}$ on $C_1$, $D'$ is well defined
and is the desired developing map. \qed{weilopen}

Using the developing map, we can pull back the geometric structure on
$X$ to a geometric structure on $\tilde{M}$. Since $\dev$ is
equivariant this geometric 
structure will descend to a geometric structure on $M$. We will refer
to this structure as $M_\dev$ to distinguish it from the
differentiable manifold, $M$. We also define the diffeomorphism,
$f_\dev: M \longrightarrow M_\dev$. The pair, $(M_\dev, f_\dev)$, is
{\em marked} geometric structure on $M$. (Often the developing map,
$D$, will be indexed by some $i$. In this case the associated
structure is $M_i$ and the diffeomorphism is $f_i$.)
With this definition one define of marked geometric structures where,
$(M_1, f_1) \sim (M_2, f_2)$ if there 
exists a geometry preserving homeomorphism, $f: M_1 \longrightarrow
M_2$, such that $f$ is isotopic to $f_2 \circ f^{-1}_1$. The space of
equivalence classes is the {\em Teichm\"uller space} of
$(G,X)$-structures on $M$. This notion of Teichm\"uller space is a
generalization of the Teichm\"uller space of complex structures
discussed in \S \ref{teichspace}.

Given a diffeomorphism, $f$, between two Riemannian manifolds, $M_0$
and $M_1$ we let $K(f)$ be the minimal biLipshitz constant of $f$. In
general $K(f)$ may be infinite, although for compact manifolds it will
always be finite. If there exists an $f$ such that $K(f)$ is finite we
say that $M_0$ and $M_1$ are {\em quasi-isometric}.

We next define geometric convergence. This concept is important
because there are sequences of manifolds that are not quasi-isometric,
yet are geometrically very close. Let $(M_n, \omega_n)$ be a sequence
of Riemannian 
manifolds with basepoints. The sequence, $(M_n, \omega_n)$, converges
to $(M_\infty, \omega_\infty)$ {\em geometrically} if there exists an
exhaustion of $M_\infty$ by compact submanifolds, $M^n_\infty$, and
$K_n$-biLipshitz maps, $f^n:M^n_\infty \longrightarrow M_n$, with
$f^n(\omega) = \omega_n$ and $K_n \rightarrow 1$.

We then have the following theorem:
\begin{prop}
\label{geomtop}
Let $D_n$ be a sequence of developing maps such that $D_n \rightarrow
D_\infty$ in $\cD(M)$. Let $\omega_n = f_{n}(\omega)$ where $\omega
\in M$. Then $(M_{n}, \omega_n) \rightarrow (M_{\infty},
\omega_\infty)$ geometrically.
\end{prop}

{\bf Proof.} First we assume that $M$ is compact. Let $V$ be a
neighborhood of $x \in \tilde{M}$ such that $D_\infty$ restricted to
$V$ is a homeomorphism onto its image. Then for large $n$, $D_n$
restricted to $V$ will be a homeomorphism. Since $D_n$ converges to
$D_\infty$ in the $C^\infty$-topology, $D_n \circ D^{-1}_\infty$
restricted to $D_\infty(V)$ converges to the identity in the
$C_\infty$ topology. Therefore $K(D_n \circ
D^{-1}_\infty)$ restricted to $D_\infty(V)$ converges to
one. Compactness then implies that $K(f^n)$ converges to one where
$f^n = f_n \circ f^{-1}_\infty$.

Let $M^n$ be an exhaustion of $M_\infty$ by compact sets
with $\omega_\infty \in M^n$ for all $n$. The biLipshitz constants of $f^n$
restricted to $M^n$ may not limit to one.  We need to relabel $f^n$ and
$M^n$ for this to happen. Using the result for compact manifolds, we
can define an strictly increasing function, $L: \integers^+
\longrightarrow \integers^+$, such that $f^n$ restricted to
$M^m$ is a $(1+1/m)$-biLipshitz map if $n > L(m)$. Now let $g^n = f^{n
+ L(1)}$ and let $X^n = M^m$ where $L(m)< n + L(1) \leq L(m+1)$. Then
$X^n$ exhaust $M_\infty$, while the $g^n$ restricted to $X^n$ have
biLipshitz constant approaching one.
\qed{geomtop} 

\section{Torus ends}
\label{torus}
The main goal of this paper is to understand the space of hyperbolic
structures on an open, tame 3-manifold $M$. To make this space
tractable we need to restrict the possible geometries of the ends. In
the next two sections we describe the allowable geometries of the ends.

Recall that an end, $E$, is a 3-manifold with boundary homeomorphic to
$S \times [0, \infty)$ where $S$ is a closed surface. As we will only
be interested in the geometry of an end ``near infinity'' we make the
following definition: Hyperbolic structures, $E_1$ and $E_2$, on $E$
are {\em isometric ends} if there exists a homeomorphism, $\phi:E_1
\longrightarrow E_2$, that is an isometry outside of a compact subset
of $E_1$. By Theorem \ref{weilopen} if $E'_1$ is a small deformation
of $E_1$ then there exists a small deformation $E'_2$ of $E_2$ such
that $E'_1$ and $E'_2$ are isometric ends.

In this section we restrict to the case where $E = T \times
[0,\infty)$ where $T$ is a torus. Let $U = \{z \in \cx | \Im z > 
0\}$ be the upper half plane of $\cx$. For each $(a,b) \in \cx
\times U$ we define a $(PSL_2\cx, \hthree)$-developing map,
$D_{(a,b)}$, in $\cD(E)$. The hyperbolic structures 
defined by these developing maps are {\em cone structures} on $E$. 

We parameterize the universal cover of $E$ as
$\reals^2 \times [1,\infty)$. Let $\gamma_1$ and $\gamma_2$ be
generators of $\pi_1(E) = \integers^2$. Then the
action of $\pi_1(\torus)$ on the universal cover is defined
by $\gamma_1(x,y,t) = (x+1,y,t)$ and $\gamma_2(x,y,t) = (x,y+1,t)$ for
all $(x,y,t) \in \reals^2 \times [1, \infty)$.

Let $s=(a,b)$ with $a \neq 0$. Before we describe $D_s$ we
define maps $\phi_s : \reals^2 \longrightarrow \cx$ by 
$$ \phi_s(x,y) = -z_0e^{xa + yab} $$
where $z_0 = \frac{1}{1-e^a}.$ Then
$$\dev_s(x,y,t) = \frac{|\phi_s(x,y)|}{\sqrt{t^2+|\phi_s(x,y)|^2}}
\left(\phi_s(x,y), t\right) + (z_0,0).$$
The corresponding holonomy representation,
$\rho_s$, is defined on the generators, $\gamma_1$ and $\gamma_2$, by
$\rho_s(\gamma_1) = e^az+1$ and $\rho_s(\gamma_2) = e^{ab}z +
\frac{e^{ab}-1}{e^a-1}$. 

Now let $a=0$ and define
$$\dev_s(x,y,t) = (x + b y, t)$$
with holonomy representation $\rho_s(\gamma_1) = z + 1$ and
$\rho_s(\gamma_2) = z + b$. 

If $s_n = (a_n, b_n)$ with $a_n \arrow 0$ and $b_n \arrow b$ as $n
\arrow 0$ then the 
maps $\dev_{s_n}$ converge to $\dev_s$ in the $C^\infty$ topology and
the holonomy representations $\rho_{s_n}$ converge to $\rho_s$.

We define $E_s$ to be the corresponding hyperbolic cone-manifold and
let $CT(E)$ be the space of ends, $E_s$. The topology
of $\cD(E)$ induces a topology on $CT(E)$.

\begin{prop}
\label{ctpar}
The space $CT(E)$ is parameterized by $\cx/ \{\pm 1\} \times
\upper$.
\end{prop}

{\bf Proof.} If two hyperbolic structures, $E_s$ and $E_{s'}$, are
isometric then they must have conjugate holonomy
representations. We separate these holonomy representations into two
classes, parabolic representations and hyperbolic representations. The
parabolic case is the easiest for $\rho_s$ will be parabolic if and
only if $s = (0,b)$ and this holonomy will be conjugate to $\rho_{s'}$
if and only if $s' = (0,b)$.

If $\rho_s$ and $\rho_{s'}$ are hyperbolic then it is possible for
$\rho_s$ and $\rho_{s'}$ to be conjugate and $E_s$ and $E_{s'}$ to not
be isometric ends. To distinguish between $E_s$ and $E_{s'}$ we
examine the intrinsic geometry of the ends. A hyperbolic holonomy
representation will fix an axis in $\hthree$. Let $T_\epsilon$ be the
torus in $E_s$ whose image under the developing map is the locus of
points distance $\epsilon$ from the axis of $\rho_s$. Similarly define
$T'_\epsilon$ in $E_{s'}$. The ends $E_s$ and $E_{s'}$ will be
isometric if and only if there is an isometry between $T_{\epsilon}$
and $T'_{\epsilon}$ in the correct homotopy class. In particular, the
length of $\gamma_1$ and $\gamma_2$ should be the same on both
$T_\epsilon$ and $T'_{\epsilon}$. These lengths can be easily
calculated using hyperbolic trigonometric functions and it is seen that
this will be the case if and only if $s=(a,b)$ and $s'=(\pm
a,b)$. \qed{ctpar}

We can now define a complex length functions, $$\cL_{x,y}:CT(\torus)
\longrightarrow \cx/\{\pm 1 \},$$ for every $x,y \in \reals$ by 
$\cL_{x,y}(E_s) = \pm a(x + by)$. By Proposition \ref{ctpar},
$\cL_{x,y}$, is well defined. If $(x,y) = (p,q)$, for relatively
prime $p$ and $q$, then $\cL_{p,q}(E_s)$ will be the length and
rotation of the holonomy representation of the $(p,q)$ curve in $E_s$.

Let $\cR_0(T) \subset \cR(T)$ be the subvariety of $\cR(T)$ consisting
of representations with infinite image and $R_0(T)$ be the quotient,
$\cR_0(T)/PSL_2\cx$.

\begin{prop}
\label{torusonto}
The holonomy map, $h:CT(E) \longrightarrow R_0(T)$, is a complex analytic
local homeomorphism.
\end{prop}

{\bf Proof.} It is well known that $R_0(T)$ is a 2-dimensional complex
manifold. By our explicit description of $\rho_s$ we see that $h$ is
holomorphic and locally injective. Since $CT(E)$ and $R_0(T)$ have the
same dimension $h$ must be a local homeomorphism. \qed{torusonto}

Proposition \ref{torusonto} allows to us define length
functions, $\cL'_{x,y}$, in neighborhood of a representation $\rho \in
R(\torus)$ by choosing a lift of this neighborhood in
$CT(\torus)$. However, in this case the $\cL'_{x,y}$ will not be
unique and instead depends on the choice of lift. In practice, the
geometric structure defining $\rho$ will determine what lift to take.

The complex length functions determine a map $$CT(\torus) \longrightarrow
\reals^2 / \{ \pm 1\} \cup \infty.$$ If $s=(a,b)$ and $a \neq 0$
then there is a unique $\pm (x,y) \in \reals^2 / \{ \pm 1\}$ such that
$\cL_{x,y}(E_s) = \pm 2 \pi \imath$. Then $\pm (x,y)$ are the {\em
hyperbolic Dehn filling coordinates} for $E_s$. If $a = 0$ the Dehn
filling coordinate for $E_s$ is $\infty$.

Let $\bar{E}_s$ be the metric completion of $E_s$. The Dehn filling
coordinate determines the topology of $\bar{E}_s$. There are three
cases.

{\bf Case 1.} If the Dehn filling coordinate for $E_s$ is $\infty$,
then $E_s$ is already complete so $\bar{E}_s =
E_s$. We refer to $E_s$ as a  {\em cusp}.

{\bf Case 2.} If the Dehn filling coordinates are $(x,y)$ and there
exists a $\theta \in \reals^+$ such that $\frac{2 \pi}{\theta} (x,y) =
(p,q)$ where $p$ and $q$ are relatively prime integers, then
$\bar{E}_s$ is an open solid torus. The metric along the core curve of
the solid torus is singular and the cross section of this curve is a
hyperbolic cone with {\em cone angle}, $\theta$. Then the curve
$(p,q)$ is the {\em meridian} of $E_s$ and is trivial in the solid
torus. In this case, we call $E_s$ a {\em rational} cone-manifold. If
the cone angle of $E_s$ is $2 \pi$ then metric extends to a smooth
metric on the solid torus $\bar{E}_s$. 

{\bf Case 3.} If $x/y$ is irrational then $E_s$ is completed by the
addition of a single point. In this case $\bar{E}_s$ is not a manifold
for the boundary of an $\epsilon$-neighborhood of the additional point
is a torus. In this case $E_s$ is an {\em irrational cone-manifold}.

\begin{lemma} 
\label{continuity}
Let $Z$ be a 1-dimensional complex submanifold of $CT(E)$ with
$E_s \in Z$. If $E_s$ is a cone-manifold with Dehn filling coordinates
$(x_0, y_0)$, then assume $\cL_{x_0,y_0}|_Z$ is injective at $E_s$. If
$E_s$ is a cusp assume there is some $(x_0,y_0)$ such that
$\cL_{x_0,y_0}|_Z$ is injective at $M_s$. In either case, $Z$ is
parameterized by Dehn filling coordinates at $E_s$.
\end{lemma}

{\bf Proof.} We first show that the map to Dehn filling coordinates is
continuous. This is clearly true away from $\infty$. Assume the sequence,
$\{E_{s_i}\}$, converges to a cusp, $E_s$, with $s_i = (\pm a_i,
b_i)$. Let $(x_i,y_i)$ be the Dehn filling coordinates for
$E_{s_i}$. By definition $a_i(x_i + b_i y_i) = \pm 2\pi \imath$. Since
$a_i \rightarrow 0$, we must have $(x_i + b_i y_i) \rightarrow
\infty$. Furthermore, the $b_i$ are bounded so $x_i$ or $y_i$ must go
to infinity, proving continuity.

We now finish the proof if $E_s$ is a cusp.  Let $A$ and $B$ be
coordinate functions for  $\bar{Z}$,  the lift of $Z$ to $\cx \times
U$ with  $(A(0), B(0)) = s$. The complex length functions,
$\cL_{x,y}$, lift to functions, $\bar{\cL}_{x_0,y_0}:(\cx \times U)
\longrightarrow \cx$ with $\bar{\cL}_{x_0,y_0}|_{\bar{Z}}$ injective
at $s$. The derivative of $\bar{\cL}_{x,y}(A(z), B(z))$ at $0$ is
$A'(0)(x + y B(0))$. By assumption, when $(x,y) = (x_0, y_0)$ this
derivative is non-zero so $A'(0) \neq 0$. If $(x,y) \neq (0,0)$ then
$(x + y B(0)) \neq 0$ because $\Im B(0) > 0$. Therefore
$\bar{\cL}_{x,y}|_{\bar{Z}}$ is injective at $(A(0), B(0))$ for all
$(x,y) \neq (0,0)$.

We now find a neighborhood, $\bar{V}$, of $(A(0), B(0))$ in $\bar{Z}$
such that $\bar{\cL}_{x,y}$ is injective on $\bar{V}$ for all
$(x,y)$. From the compactness of the unit circle we can first find
such a $\bar{V}$ for all $(x,y)$ on the unit circle. We then note that
if $\bar{\cL}_{x,y}$ is injective on $\bar{V}$ then $\bar{\cL}_{kx,
ky}$ is also injective on $\bar{V}$ if $k \neq 0$ so the same
neighborhood, $\bar{V}$, works for all $(x,y)$. 

If necessary we shrink $\bar{V}$ so that it is a branched double cover
of a neighborhood $V$ in $Z$. The injectivity of the $\bar{\cL}_{x,y}$
on $\bar{V}$ imply that the $\cL_{x,y}$ are injective on $V$.
Therefore $\cL_{x,y}(E_t) = \pm 2 \pi \imath$ for at most one $E_t$ in
$V$ and $\cL_{x,y}$ is zero only at $E_s$. This implies that the map
to Dehn filling coefficients is injective. Since it is also continuous
the lemma is proved in the cusped case.  

The proof is essentially the same if $E_s$ is a cone-manifold except
the injectivity of $\bar{\cL}_{x_0,y_0}$ only implies that length
functions with coefficients near $(x_0,y_0)$ are injective. However,
this is enough to prove the lemma. See \cite{Hodgson:Kerckhoff:cone},
Theorem 4.8, for a slicker version of this proof.
\qed{continuity}

\section{Higher genus ends}
\label{ends}
\subsection{Geometrically finite ends}
In \S \ref{findef} we defined the notion of geometric finiteness
for ends of complete hyperbolic manifolds. We now generalize this
definition so that it will apply to ends of cone-manifolds. Let $E = S
\times [0,1)$. We give the union, $E \cup S$ the product topology, $S
\times [0,1]$. A developing map, $D$, for $E$ is
{\em geometrically finite without rank one cusps} if $D$ extends to a
local homeomorphism, $\bar{D}:(\tilde{E} \cup \tilde{S})
\longrightarrow (\hthree \cup \chat)$. The restriction of $\bar{D}$ to
$S \times \{1\}$ will be a $(PSL_2\cx, \chat)$-developing map for
$S$. A $(PSL_2\cx, \chat)$-structure on $S$ is a {\em projective
structure}. If $D$ is globally injective, the holonomy of $D$
will be discrete and faithful and this definition of geometric
finiteness agrees with that given in \S 2. We also remark that a
projective structure on $S$ determines a conformal structure. The
restriction of $\bar{D}$ to $S$ defines a {\em projective structure at
infinity} for $E$ and the underlying conformal structure is the
conformal structure at infinity defined in \S \ref{findef}.
The following proposition is proved in \cite{Bromberg:rigidity}.

\begin{prop}
\label{gfisom}
To hyperbolic structures $E_1$ and $E_2$ on $E$ are isometric ends if
and only if they have the same projective structure at infinity.
\end{prop}

Before we continue our discussion of geometrically finite ends we need
to review 
some of the basic facts about projective structures on $S$. We denote
the Teichm\"uller space of projective structures, $P(S)$. As mentioned
above, a projective structure determines a complex structure, so there
is a projection, $P(S) \longrightarrow T(S)$, where $T(S)$ is the
Teichm\"uller space of complex structures discussed in \S
\ref{teichspace}.  The fibers of this projection are the holomorphic
quadratic differentials on each
complex structure and therefore $P(S)$ has complex dimension
$6g-6$. We also have the following fundamental theorem: 
\begin{theorem}[Hejhal]
\label{hejhal}
The holonomy map $P(S) \longrightarrow R(S)$ is a complex, local
 homeomorphism.
\end{theorem}

We can realize a conformal structure, $\mu$, on $S$ as the quotient,
$U/\Gamma$, where $\Gamma$ is a Fuchsian group. If we use $U$ as our
topological model for $\tilde{S}$ with $\Gamma$ the group of deck
transformations then any projective structure on $S$ with conformal
structure structure $\mu$ is determined by a conformal developing map
$f: U \longrightarrow \chat$. In this way the pair $(f, \Gamma)$
determines a projective structure on $S$. Note that a projective
structure does not determine a unique pair, $(f, \Gamma)$.

{\bf Remark.} If $(f', \Gamma')$ is another projective structure we
cannot compare $f$ to $f'$ because although they have the same domain,
$U$, this domain is acted on by a different group of deck
transformations for each map. On the other hand if $(f, \Gamma)$ and
$(f', \Gamma')$ 
are near each other in $P(S)$ then there will be a diffeomorphism,
$\phi: U \longrightarrow U$, such that $f' \circ \phi$ is a developing
map with deck transformations $\Gamma$ and with $f$ and $f' \circ
\phi$ close in the $C^\infty$ topology.

\medskip

We now construct a geometrically finite end, $E_{(f, \Gamma)}$, that has
projective structure at infinity $(f, \Gamma)$. To do so we will use a
general method for extending a conformal map, $f:\upper
\longrightarrow \chat$, to a map $\Theta_f: \hthree \longrightarrow
\hthree$. For each $z \in \upper$ we let $M^f_z$ be the {\em 
osculating} Mobius transformation, i.e. the unique element of
$PSL_2\cx$ whose two-jet agrees with $f$ at $z$. Let $P \subset
\hthree$ be the hyperbolic plane bounded by $\reals \cup
\infty$. Orient $P$ such that the normal projection onto $\upper$ is
orientation preserving and let $P_t$ be the image of the time $t$
normal flow of $P$. For every $x \in \hthree$ there is a unique
geodesic, $r$, orthogonal to $P$ through $x$. Let $r(x)$ be the limit
of this geodesic in $\upper$ and define $\Theta_f(x) =
M^f_{r(x)}(x)$.

Unfortunately $\Theta_f$ will not be a local homeomorphism on
all of $\hthree$. To examine when $\Theta_f$ is a local homeomorphism
we will need a result of Epstein and Anderson. To state it we need the
notion of the {\em Schwarzian derivative} of a locally injective
holomorphic function,
$$ SCf(z) = \left(\frac{f''}{f'}\right)' -
\frac12\left(\frac{f''}{f'}\right)^2, $$
and its hyperbolic norm in $\upper$,
$$ \| SCf(z) \| = z^2|SCf(z)|.$$
If $f$ is conformal developing map for $S$, then a
straightforward calculation shows that $\| SCf(z) \| = \|
SCf(\gamma(z)) \|$ for all $\gamma \in \pi_1(S)$.

\begin{theorem}
\label{epstein}
Let $p \in P_d$, and $e_3$ the oriented unit normal vector for $P_d$
at $p$. Then there exists a orthonormal basis $e_1$ and $e_2$ for the
tangent space of $P_d$ at $p$ and an orthonormal basis $e'_1$, $e'_2$
and $e'_3$ for the tangent space of $\hthree$ at $\Theta_f(p)$
such that the derivative of $\Theta_f$ at $p$ in terms of these two bases is 
$$\left( \begin{array}{ccc}
1 + \frac{\|SCf(r(p))\|}{\cosh d} & 0 & 0 \\
0 & 1 - \frac{\|SCf(r(p))\|}{\cosh d} & 0 \\
0 & 0 & 1 \\
\end{array}
\right)_\cdot$$
\end{theorem}

{\bf Remark 1.} In \cite{Anderson:projective} what is actually
calculated is the derivative of 
the composition of the map from $\upper$ to $P_d$ and
$\Theta_f$. However, since the map from $\upper$ to $P_d$ just
multiplies the metric by $\cosh d$ our result follows immediately from
the work in \cite{Anderson:projective}.

{\bf Remark 2.} The Schwarzian derivative is a holomorphic quadratic
differential. A quadratic differential defines a pair of singular
foliations, the horizontal and vertical trajectories. The $r_* e_1$
will lie in the direction of the horizontal trajectory while $r_* e_2$
will lie in the direction of the vertical trajectory. The foliations
will be singular were the Schwarzian is zero. However, in this case
the choice of $e_1$ and $e_2$ is arbitrary.
\medskip

We have the following corollary of Epstein and Anderson's result:

\begin{cor}
\label{thetahomeo}
If $f$ is equivariant, as above there exists a $d_0$ such that if $d >
d_0$ and $x \in P_d$ then $\Theta_f$ is a local homeomorphism at $x$.
\end{cor}

{\bf Proof.} We first note that $\Theta_f$ is equivariant. This
follows immediately from the relationship $M^f_{r(\gamma(x))} =
\rho_s(\gamma) \circ M^f_{r(x)} \circ \gamma^{-1}$ for all $x \in
\hthree$ and $\gamma \in \pi_1(S)$. Furthermore since $S$ is
compact, equivariance also implies that $\|SCf(z)\|$ is bounded on
$\upper$. We then choose $d_0$ such that $\cosh d_0 = \sup_{z \in
\upper} \|SCf(z)\|$ proving the corollary. 
\qed{thetahomeo}

Now choose $d_0$ as in Corollary \ref{thetahomeo} and let $H_{d_0}$ be
the the subset of hyperbolic space bounded by $P_{d_0}$ and
$\upper$. Then the universal cover, $\tilde{E}$ can be identified with
$H_{d_0}$ and the deck transformations acting on $H_{d_0}$ will be
$\Gamma$. Then $\Theta_f$ 
restricted to $H_{d_0}$ is a developing map for $E$ that defines a
geometrically finite structure $E_{(f,\Gamma)}$ on $E$ with
projective structure at infinity $(f, \Gamma)$.

The following proposition shows that two geometrically finite ends are
close if there projective structures at infinity are close.

\begin{prop}
\label{quasiends}
Let $E_{(f_i,\Gamma_i)}$ be a sequence of geometrically finite ends
such that $(f_i, \Gamma_i)$ converges to $(f_\infty,\Gamma_\infty)$ in
$P(S)$. Then there exists $K_i$-biLipshitz homeomorphisms,
$f_i:E_{(f_\infty,\Gamma_\infty)} \longrightarrow E_{(f_i,\Gamma_i)}$,
with $K_i \rightarrow 1$.
\end{prop}

{\bf Proof.} 
If $(f_i,\Gamma_i) \rightarrow (f_\infty,\Gamma_\infty)$
in $P(S)$  then we can choose the representatives for $(f_i,
\Gamma_i)$ such that $f_i \rightarrow f_\infty$ in the $C^\infty$
topology and such that $\Gamma_i \rightarrow \Gamma_\infty$ in
$\cR(S)$. For each $f_i$ let $d_i$ be the constant given Corollary
\ref{thetahomeo}. Since $f_i \rightarrow f_\infty$, $d_i \rightarrow
d_\infty$. Choose $d'$ such that $d' > d_i$ for large $i$. To prove
the proposition we construct developing maps, $D_i:H_{d'}
\longrightarrow \hthree$, with the fixed group of deck transformations,
$\Gamma_\infty$, and show that $(D_\infty)_* \circ (D_{i}^{-1})_*$
converges uniformly to the identity in the hyperbolic metric.

Since $\Gamma_i \rightarrow \Gamma_\infty$ there exist conjugating
diffeomorphisms $\phi_i:\upper \longrightarrow \upper$ such that
$\Gamma_i = \phi_{i}^{-1} \circ \Gamma_\infty \circ \phi_i$ and $\phi
\rightarrow \id$. We can extend the $\phi_i$ to maps $\Phi_i:\hthree
\longrightarrow \hthree$. The map $\Phi_i$ is determined uniquely by
the two conditions that $\Phi_i$ takes $P_d$ to itself and
$\phi_i(r(x)) = r(\Phi_i(x))$. For the hyperbolic metric on $\upper$,
$\phi$, is $K'_i$-biLipshitz with $K'_i \rightarrow 0$. Since $P_d$ is
a scaled version of the hyperbolic plane, $\Phi_i$, will also be
$K'_i$-biLipshitz and the derivatives, $(\Phi_i)_*$ will converge
uniformly to the identity.

For the sequence of conformal map, $f_i$, the derivative of
$\Theta_{f_i}$ will converge uniformly along a ray perpendicular to
$P$ by Theorem \ref{epstein}. Choose a fundamental domain for the
action of $\Gamma_\infty$ on $\upper$ and let $F$ be the pre-image of
this domain under $r$. Then $F$ is a compact family of rays
perpendicular to $P$ so $(\Theta_{f_i})_*$ will converge uniformly to
$(\Theta_{f_\infty})_*$ on $F$.

Now define $D_\infty = \Theta_{f_\infty}|_{H_{d'}}$ and $D_i =
\Theta_{f_i} \circ \Phi_i|_{H_{d'}}$. Then $(D_\infty)_* \circ
(D_{i}^{-1})_*$ will converge uniformly on $F$ and therefore by
equivariance on all of $H_{d'}$. \qed{quasiends}

We let $GF(E)$ denote the space of structures $E_{(f,\Gamma)}$.
We have shown that every projective structure is the projective
structure of a geometrically finite end. This fact along with
Propositions \ref{gfisom} and \ref{quasiends} and Theorem \ref{hejhal}
implies the following proposition:
\begin{prop}
\label{gfonto}
The holonomy map $GF(E) \longrightarrow R(S)$ is locally onto.
\end{prop}

We close this section by defining a map, $b:GF(E)
\longrightarrow T(S)$ by $b(E_{(f, \Gamma)}) = \mu$ where $\mu$ is the
conformal structure of $\upper/\Gamma$. If $V$ is a neighborhood of
the holonomy representation, $\rho$, of $E_{(f,\Gamma)}$ such that $h|_V$ is
injective then Theorem \ref{hejhal} allows us to define a map, $b':V
\longrightarrow T(S)$ by composing $b$ with the local inverse of $h$.
Since $\rho$ will not be the holonomy of a unique structure in
$GF(E)$, the map $b'$ will depend on both $\rho$ and $E_{(f,\Gamma)}$.

\subsection{Geometrically infinite ends}
If the hyperbolic structure on $E$ is not geometrically finite, then
it is not know in general if $E$ has local deformations. We will
examine a special case. If the developing map, $D$, for $E$ is
injective then the image of the holonomy of $E$ will be a Kleinian
group $\Gamma$. The quotient, $\hthree/\Gamma = M$, will be a
hyperbolic 3-manifold homeomorphic to $S \times \reals$ and $E$ will
isometrically embed as an end of $M$. Following our previous notation,
we let $QD(\Gamma; S)$ be the space of quasiconformal deformations of
$\Gamma$ that fix the end-invariant of $E$. Now for each Kleinian
group $\Gamma' \in QD(\Gamma; S)$, the quotient $\hthree/\Gamma' = M'$
will again be a hyperbolic 
structure on $S \times \reals$. However, it is not immediate that the
ends of $M'$ are near the ends of $M$ when $\Gamma'$ is near
$\Gamma$. To show this we need to extend a quasiconformal map of
$\chat$ to a biLipshitz map of $\hthree$. The following result can be
found in \cite{McMullen:book:RTM}.

\begin{theorem}
\label{extmap}
Let $\phi: \chat \longrightarrow \chat$ be a $k$-quasiconformal
homeomorphism. Then $\phi$ extends continuously to a $K$-biLipshitz
map, $\Phi: \hthree \longrightarrow \hthree$, where $K =
k^{3/2}$. Furthermore, this extension is natural in the sense that if
$\gamma$ and $\gamma'$ are elements of $PSL_2\cx$ such that
$$\phi \circ \gamma' = \gamma \circ \phi$$
then
$$\Phi \circ \gamma' = \gamma \circ \Phi.$$
\end{theorem}

Given $\Gamma' \in QD(\Gamma; S)$ there is a
$k$-quasiconformal map, $\phi:\chat \longrightarrow \chat$, 
that conjugates $\Gamma$ to $\Gamma'$. Let $\Phi$ be the
extension of $\phi$ given by Theorem \ref{extmap}. If we identify the
universal cover, $\tilde{E}$, of $E$ with a subspace of $\hthree$ such
that $\Gamma$ acts as deck transformations then $\Phi|_{\tilde{E}}$ is
a developing map that defines an end $E'$ near $E$. Let $QD(E)$ be the
space of ends constructed through such developing maps with the
developing map topology. Theorem \ref{extmap} implies that these
developing maps will vary continuously in the space of all developing
maps so we have:

\begin{prop}
\label{infonto}
The holonomy map, $QD(E) \longrightarrow QD(\Gamma; S)$, is a
homeomorphism.
\end{prop}

\section{The Dehn filling theorems}
Let $M$ be the interior of a compact, hyperbolizable 3-manifold
$N$. Assume $\del N$ has $n$ components consisting of tori, $T_1,
\dots, T_k$, and higher genus surfaces, $S_{k+1}, \dots, S_n$. Given a
compact core $M_0$ of $M$ we label the components of $M \backslash
M_0$, $E_1, \dots, E_n$. If $i \leq k$, then $E_i$ will be homeomorphic
to $T_i \times [0,1)$ and if $i > k$ then $E_i$ will be homeomorphic to
$S_i \times [0,1)$ where $S_i$ has genus $> 1$.

A hyperbolic structure on $M$ is a {\em cone-manifold} if for $i \leq
k$ there exists an $E_s \in CT(E_i)$ such that $E_i$ and $E_s$ are end
isometric. $M$ is a {\em rational} cone-manifold if all of these
structures are rational or cusps. $M$ is {\em geometrically finite
without rank one cusps} if for $i>k$ there exists $E_{(f, \Gamma)} \in
GF(E_i)$ such that $E_i$ and $E_{(f, \Gamma)}$ are end isometric.

Let $\cG\cF(M)$ be the space of developing maps that induce
geometrically finite cone structures and let $GF(M)$ be the associated
Teichm\"uller space. The topology on $\cG\cF(M)$ induces a topology on
$GF(M)$. To see this topology more explicitly we define a sub-basis
for it. Let $M' \in GF(M)$ be a geometrically finite cone structure
and let $M'_0$ be a compact core for $M'$. Then $M'(\epsilon)$ is the
set of structures in $GF(M)$ such that there exists an $(1 +
\epsilon)$-biLipshitz embedding of $M'_0$ in the correct homotopy
class. The sets $M'(\epsilon)$ form a sub-basis of neighborhoods for
$M'$.

Now assume $M$ is a rational cone manifold with holonomy
representation $\rho$. If $E_i$, with $i \leq k$, is not a cusp, let
$(p_i, q_i)$ be the meridian. Otherwise choose $(p_i,q_i)$ to be any
simple closed curve on $T_i$. We then have length functions,
$\cL_{p_i,q_i}:CT(E_i) \longrightarrow \cx/\{\pm 1\}$, and
$\cL_{p_i,q_i}':V_i \longrightarrow \cx/\{ \pm1\}$, where $V_i$ is a
neighborhood of the restriction of $\rho$ to $E_i$.

If $i>k$ we have maps $b_i:GF(E_i) \longrightarrow T(S_i)$ and
$b'_i:V_i \longrightarrow T(S_{g_i})$ where $V_i$ is a neighborhood of
the restriction of $\rho$ to $E_i$. We
combine these maps in the obvious way to get two maps, $\Phi:GF(M)
\longrightarrow \cx^k \times T(S_{k+1}) \times \cdots \times
T(S_n)$, and $\Phi':V \longrightarrow \cx^k \times
T(S_{g_{k+1}}) \times \cdots \times T(S_{g_n})$ where $V$ is a
neighborhood of $\rho$. The main result of \cite{Bromberg:rigidity} is:

\begin{theorem}
\label{rephomeo}
If $M$ is a geometrically finite rational cone-manifold with holonomy
representation, $\rho$, and all cone angles $\leq 2\pi$, then $\Phi'$
is a local homeomorphism at $\rho$. 
\end{theorem}

We now prove:

\begin{theorem}
\label{geomhomeo}
If $M$ is a geometrically finite rational cone-manifold with all cone
angles $\leq 2\pi$, then $\Phi$ is a local homeomorphism at
$M$. Therefore $GF(M)$ is locally parameterized by hyperbolic
Dehn filling coordinates and the conformal structure at infinity.
\end{theorem}

{\bf Proof.} Let $h:GF(M) \longrightarrow R(M)$ be the holonomy
map. By the definition of $\Phi'$, $\Phi = \Phi' \circ h$. To
complete the proof we need to show that $h$ is a local homeomorphism.

Let $D \in GF(M)$ be a developing map for $M$, $\rho = h(M)$ the
holonomy and $M_0$ a compact core. If $\rho' \in R(M)$ is near $\rho$
then by Theorem \ref{weils theorem} there exist $D' \in \cD(M_0)$ such
that the holonomy of $D'$ is $\rho'$ and $D'$ is near $D|_{M_0}$. By
Theorem \ref{weilopen} $D'$ extends to a developing map in $\cG\cF(M)$
since by Proposition \ref{torusonto} and Proposition \ref{gfonto} for
all nearby representations there is a always a cone on $E_i$ if $i
\leq k$ and a geometrically finite structure without rank one cusps if
$i > k$. Therefore $h$ is onto.

Now assume that $M^1$ and $M^2$ are near $M$ and that $h(M^1) =
h(M^2)$. We again apply Theorem \ref{weils theorem} to find a compact
core $M^{1}_0$ of $M^1$ that isometrically embeds in $M^2$. Since
$M^{1}_0$ will extend to a unique geometrically finite cone structure
this implies that $M^1 = M^2$ and $h$ is injective.

Therefore $h$ is a local homeomorphism and combining this fact with
Theorem \ref{rephomeo} implies that $\Phi$ is a local
homeomorphism. 

Lemma \ref{continuity} allows us to replace the complex lengths with
Dehn filling coefficients in the parameterization.
\qed{geomhomeo}

We can now prove Thurston's Dehn filling theorem, generalized to
geometrically finite manifolds without 
rank one cusps. Before we state it we recall the notion of topological
Dehn filling. 

Let $M$ be a compact 3-manifold whose boundary contains a torus
$T$. Choose a basis for $\pi_1(T) = \integers \oplus \integers$ so
that an element of $\pi_1(T)$ is determined by a pair of integers
$(p,q)$. For each relatively prime pair, $(p,q)$, there is a manifold,
$N(p,q)$ and an embedding, $f_{p,q}:M \longrightarrow M(p,q)$, such
that $M(p,q) - f_{p,q}(M)$ is a solid torus with boundary $f_{p,q}(T)$
and the image of the $(p,q)$ curve in $N(p,q)$ is trivial. $N(p,q)$ is
the {\em $(p,q)$-Dehn filling} of $M$. Let $\gamma$ denote the core
curve of the solid torus. If $M(p,q)$ has a hyperbolic
structure on its interior, then $M(p,q)$ is a {\em hyperbolic Dehn
filling} of $M$ if $f_{p,q}$ can be chosen such that $M(p,q) -
f_{p,q}(M)$ is the geodesic representative of $\gamma$. Note that a
hyperbolic structure on $M(p,q)$ may not be a hyperbolic Dehn filling
of $M$ if the geodesic representative of $\gamma$ is not isotopic to
$\gamma$. Also note that the holonomy representation $\rho$ of
$M(p,q)$ also defines a non-faithful, holonomy representation,
$\hat{\rho}$, for $M$ via the map $f_{p,q}$. We further note that if
$M$ also has a hyperbolic structure the map $f_{p,q}$ allows us to
compare the end-invariant of $M$ to that of $M(p,q)$.

If $M$ has $k$ torus boundary components, we can Dehn fill each of
them. Let relatively prime integers, $(p_i,q_i)$, be the Dehn filling
coefficients for the $i$-th torus and let $(\bpq) = (p_1,
q_1;\ldots;p_k,q_k)$. Then $M(\bpq)$ is the $(\bpq)$-Dehn filling of
$M$.

Define $|\bpq| = |p_1| + |q_1| + \cdots + |p_k| + |q_k|$.

\begin{theorem}
\label{bigdehn}
Let $M$ be a compact 3-manifold with $k$ torus boundary components and
assume $M$ has a geometrically finite hyperbolic
structure without rank one cusps with holonomy $\rho$. We then have
the following:
\begin{enumerate}
\item Except for
a finite number of 
pairs for each $i$, for each collection of relatively prime pairs
$(\bpq)$ there exist a geometrically finite hyperbolic $(\bpq)$-Dehn
filling $M(\bpq)$ of $M$ with $M(\bpq)$ having
the end-invariant as $M$. 

\item $\rho_{\bpq} \rightarrow \rho$ as $|\bpq| \rightarrow \infty$.

\item If $X$ is a submanifold of $M$ such that the complement of $M$
is a neighborhood of the cusps then $f_{\bpq}|_X$ is
$K_{\bpq}$-biLipshitz with $K_{\bpq} \rightarrow 1$ as $\|\bpq\|
\rightarrow \infty$.
\end{enumerate}
\end{theorem}

{\bf Proof.} (1) and (2) By Theorem \ref{geomhomeo} $GF(M)$ is locally
parameterized by the conformal structures at infinity and hyperbolic
Dehn filling coordinates. If the Dehn filling coordinates for a
structure $M' \in GF(M)$ are a set of relatively prime integers,
$(\bpq)$, then the metric completion of $M'$ is a smooth hyperbolic
structure on $M(\bpq)$. A neighborhood of $\infty \times \cdots \times
\infty$ in $(\reals^2 \cup \infty)^k$ contains all $(\bpq)$ that
exclude a finite number of pairs for each $i$.

(3) If $X$ were compact this would simply state that the
$M_{\bpq}$ converge geometrically to $M$ which is Proposition
\ref{geomtop}.  We are allowed to include the geometrically finite
ends in $X$ because of Proposition \ref{quasiends}. \qed{bigdehn}

We next prove the Dehn filling theorem for the class of geometrically
infinite 3-manifolds discussed in \S \ref{infinite}. For a manifold,
$M$, in this class, Theorem \ref{infinitedehn} describes a local
parameterization for the variety $R(M; \del_1 M)$. There is a natural
space of hyperbolic structures, $GF(M; \del_1 M)$, on $M$ associated to
the representations in $R(M; \del_1 M)$. As before, we fix the class
of structures the ends of $M$ can have. The isometry type of the torus
and geometrically finite ends is defined as above. Each geometrically
infinite end will be end isometric to an end in $QD(E_i)$.
We then have a map, $\Phi:GF(M; \del_1 M) \longrightarrow \cx^k \times
T(\del_1 M)$. Following the proofs of Theorems \ref{geomhomeo} and
\ref{bigdehn} along with Proposition \ref{infonto} we have:

\begin{theorem}
\label{infinitedehngeom}
$\Phi$ is a local homeomorphism at $M$. Furthermore (1) - (3) of
Theorem \ref{bigdehn} hold.
\end{theorem}

\section{A density theorem}
Let $AH(M) \subset R(M)$ be the space of discrete faithful
representations of $\pi_1(M)$ and $MP(M) \subset AH(M)$ those
representations whose quotient manifold is geometrically finite
without rank one cusps (often referred to as {\em minimally
parabolic}). Expanding on a more restrictive conjecture of Bers,
Thurston conjectured that $\overline{MP(M)} = AH(M)$.

If $M$ has incompressible boundary,
the density conjecture is a consequence of the ending lamination
conjecture. Thurston's relative compactness theorems imply
that every possible end invariant is realized in
$\overline{MP(M)}$. The ending lamination conjecture implies that this
structure is unique and hence $\overline{MP(M)} = AH(M)$. 

In fact the ending lamination conjecture implies that there is a
hierarchy of quasiconformal deformation spaces. Let $\Gamma$ and
$\Gamma_0$ be isomorphic Kleinian groups with ending laminations
$\lambda$ and $\lambda_0$, respectively. If $\lambda_0$ is a
sublamination of $\lambda$ then the ending lamination conjecture
implies that $QD(\Gamma) \subset \overline{QD(\Gamma_0)}$. The
particular case that we will be interested in is when $\lambda -
\lambda_0$ is a collection of isolated simple closed curves in
$\lambda$. If $\lambda$ is entirely simple closed curves then $\Gamma$
is geometrically finite and the fact that $QD(\Gamma_0) \subset
\overline{QD(\Gamma)}$ is a result of Abikoff \cite{Abikoff:bdry}. In
the following theorem we show that the isolated rank one cusps can be
resolved in geometrically infinite manifolds.

\begin{theorem}
\label{density}
Let $\Gamma$ be a Kleinian group such that the geometrically infinite
ends of $M = \hthree/\Gamma$ are incompressible and don't have any
accidental parabolics. Let $\lambda_0$  and $\lambda_1$ be the union
of the ending laminations for geometrically infinite and finite ends,
respectively. Then there exists a Kleinian group, $\Gamma_0$, with
ending lamination, $\lambda_0$, such that $QD(\Gamma) \subset
\overline{QD(\Gamma_0)}$. 
\end{theorem}

{\bf Proof.} Let $N$ be a compact topological manifold whose interior
is homeomorphic to $M$.
View $\lambda_0$ as a collection of curves in $\del N$, define $\cC$
to be a 
collection of simple closed curves in $N$ that are isotopic to the
ending lamination $\lambda_1$ and let $\hat{N} = M - \cC$. Since the
curves in $\cC$ can isotoped into $\del N$ there isotopy of $N$ inside
of itself that ``pushes'' the $\del N$ past $N$. In other words there
is a map $f: N \longrightarrow \hat{N} \subset N$ which is a
homeomorphism onto its image.

Using the Kleinian-Maskit combination theorem we can find parabolic
isometries, $\gamma_1, \dots, \gamma_n$, such that
$\hthree/\hat{\Gamma}$ is homeomorphic to $\hat{N}$ with $\hat{\Gamma}
= \Gamma*\gamma_1*\cdots*\gamma_n$. Furthermore we can choose these
parabolics such that the image of $f_*(\pi_(N))$ in $\hat{\Gamma}$ is
$\Gamma$. Since there is also an inclusion of $\hat{N}$ in $N$ we can
compare the end-invariants of $\hat{N}$ to $N$. The geometrically
infinite ends of $\hat{N}$ will be isometric to the corresponding
infinite ends of $N$.

Let $\hat{\rho}$ be the representation of $\pi_1(\hat{N})$ with image
$\hat{\Gamma}$. Choose a basis for the tori  $\del \hat{N}$ associated
to $\cC$ such that the $(1,0)$-curve is trivial in $N$. Apply Theorem
\ref{infinitedehngeom} to obtain manifolds $M_n$ that are the $({\bf
1},{\bf n})$-hyperbolic Dehn filling of $\hat{M}$ along $\cC$. Let
$f_n:\hat{M} \longrightarrow M_n$ be the embeddings obtained from the
theorem and $\hat{\rho}_n$ the representations of $\pi_1(\hat{N})$.
By this construction the $M_n$ will be homeomorphic to $M$ and the
maps $f_n \circ f$ will be homotopy equivalences. Let $\rho_n$ be the
representation of $\pi_1(N)$ induced by $f_n \circ f$. Since
$\hat{\rho}_n \rightarrow \hat{\rho}$, $\rho_n \rightarrow \rho$ since
$\rho_n$ is the restriction of $\hat{\rho}_n$ to the subgroup
$\pi_1(N)$.

Let $\Gamma_0$ be the image of $\rho_n$ for some $n$. To show that
$\Gamma \in \overline{QD(\Gamma_0)}$ we need to know that all the
$M_n$ are quasi-isometric. Since $M_n$ and $M_{n'}$ are homeomorphic,
the are quasi-isometric if the ends of $M_n$ are quasi-isometric to
the corresponding ends of $M_{n'}$. This is another consequence of
Theorem \ref{infinitedehngeom}. Therefore $\Gamma_n \in QD(\Gamma_0)$
and $\Gamma \in \overline{QD(\Gamma_0)}$.

Let $\Gamma' \in QD(\Gamma)$ and repeat the procedure. The ends of $M'
= \hthree/\Gamma'$ will be quasi-isometric to the ends of
$M$. Therefore the ends of $M'_n$ will be quasi-isometric to the ends
of $M_n$ so all the $\Gamma'_n$ will lie in $QD(\Gamma_0)$ and
therefore $QD(\Gamma) \subset \overline{QD(\Gamma_0)}$. \qed{density}

We close with a question. In Theorems \ref{infinitedehngeom} and
\ref{density} we are forced to restrict to manifolds with
geometrically infinite ends that do not have accidental
parabolics. This is because we know that these ends have
quasi-conformal deformations (Proposition \ref{bgroup}) giving the
necessary semi-stability property.  On the other hand, geometrically
finite ends are stable even if they do not have quasiconformal
deformations. Furthermore at a geometrically finite representation
there is a half dimensional subspace of $R(S)$ where the end can be
deformed and the conformal structure at infinity is fixed. By analogy,
the ending lamination is a ``degenerate'' conformal structure so we
ask if there is a half dimensional subspace of $R(S)$ where a infinite
end can be deformed preserving the quasi-isometry type.

\bibliographystyle{math}
\bibliography{math}

\begin{sc}
\noindent
Department of Mathematics\\
University of Michigan\\
East Hall, 525 E University Ave\\
Ann Arbor, MI 48109\\
\end{sc}

\end{document}